\documentclass[]{amsart}
\usepackage{amssymb}
\usepackage{mathrsfs}                     
\usepackage[all]{xy}
\allowdisplaybreaks

\def\undersetbrace#1\to#2{\underbrace{#2}_{#1}}
\def\oversetbrace#1\to#2{\overbrace{#2}^{#1}}
\def\AMSunderset#1\to#2{\underset{#1}{#2}}
\def\AMSoverset#1\to#2{\overset{#1}{#2}}

\swapnumbers

\newtheorem*{proposition*}{Proposition}

\newtheorem*{theorem*}{Theorem}
\newtheorem*{maintheorem*}{Main Theorem}

\newtheorem*{lemma*}{Lemma}

\newtheorem*{corollary*}{Corollary}

\def\ign#1{}             
\def\o{\circ}
\def\X{\mathfrak X}
\def\al{\alpha}
\def\be{\beta}

\def\th{\theta}
\def\io{\iota}

\def\ta{\tau}
\def\ph{\varphi}

\def\ps{\psi}
\def\om{\omega}
\def\Ga{\Gamma}
\def\De{\Delta}

\def\La{\Lambda}

\def\Ph{\Phi}

\def\Om{\Omega}
\def\i{^{-1}}
\def\x{\times}
\def\p{\partial}
\let\on=\operatorname
\def\L{\mathcal L}
\def\Diff{\on{Diff}}
\def\Dens{\on{Dens}}
\def\Prob{\on{Prob}}
\def\Vol{\on{Vol}}

\def\R{{\mathbb R}}

\begin{document}
\title[]
{Uniqueness of the Fisher--Rao metric on the space of smooth densities
}
\author{Martin Bauer, Martins Bruveris, Peter W. Michor}
\address{
Martin Bauer, Peter W.\ Michor: Fakult\"at f\"ur Mathematik,
Universit\"at Wien, Os\-kar-Mor\-gen\-stern-Platz 1, A-1090 Wien, Austria.\newline\indent
Martins Bruveris: 
Department of Mathematics, 
Brunel University London, Ux\-bridge, UB8 3PH, United Kingdom}
\email{bauer.martin@univie.ac.at}
\email{martins.bruveris@brunel.ac.uk}
\email{peter.michor@univie.ac.at}

\date{{\today} } 

\thanks{MB was supported by `Fonds zur
F\"orderung der wissenschaftlichen                    
Forschung, Projekt P~24625'} 
\keywords{Fisher--Rao Metric; Information Geometry;  Invariant Metrics; Space of Densities; 
Groups of Diffeomorphisms}
\subjclass[2010]{Primary 58B20, 58D15} 

\begin{abstract} 
On a closed manifold of dimension greater than one, every smooth weak Riemannian metric on the space of smooth positive probability densities, that is invariant under the action of the diffeomorphism group, is a multiple of the 
Fisher--Rao metric. 
\end{abstract}
\def\LaTeXonly{}

\maketitle

\subsection*{Introduction} 
The Fisher--Rao metric on the space $\on{Prob}(M)$ of probability densities is of importance in the 
field of information geometry. Restricted to finite-dimensional submanifolds of $\on{Prob}(M)$, 
so-called statistical manifolds, it is called Fisher's information metric \cite{Ama1985}. The 
Fisher--Rao metric has the property that it is invariant under the action of the diffeomorphism 
group. The interesting question is whether it is the unique metric possessing this invariance 
property. A uniqueness result was established \cite[p. 156]{Cen1982} for Fisher's information 
metric on finite sample spaces and \cite{AJLS2014} extended it to infinite sample spaces.

The Fisher--Rao metric on the infinite-dimensional manifold of all positive probability densities 
was studied in \cite{Fri1991}, including the computation of its curvature. 
A consequence of our main theorem in this article is the infinite-dimensional analogue of the result in \cite{Cen1982}: 

\begin{theorem*} Let $M$ be a compact manifold without boundary of dimension $\geq 2$.
Then any smooth weak Riemannian metric on the space  $\Prob(M)$ of smooth positive probability 
densities, that is invariant under the action of the diffeomorphism group of $M$, is a multiple 
of the Fisher--Rao metric.
\end{theorem*}

The situation for a 1-dimensional manifold is described at the end of the paper.
Our result holds for smooth positive probability densities on a compact manifold. However, the proof can be adapted to 
a suitable (and there are many choices) space of densities on a non-compact manifold.  
In \cite{AJLS2014} the authors prove a related result about the uniqueness of an invariant 2-tensor 
field on the space of probability densities. However they assume that the tensor is defined also on 
non-smooth densities and is invariant not only under smooth diffeomorphisms, but under all 
sufficient statistics. This is a stronger invariance assumption, allowing the authors to consider 
probability densities that are step functions, thus reducing the problem to the finite-dimensional 
case of \cite{Cen1982}.     

\subsection*{Acknowledgments}
This question was brought to our attention during a workshop at \"Oli-H\"utte above Bad Gastein in Austria, July 14--20, 2014. 
We thank all the participants of the workshop for the friendly atmosphere and helpful discussions.

\subsection*{The space of densities}
Let $M^m$ be a smooth manifold without boundary.  
Let $(U_\al,u_\al)$ be a smooth atlas for it. The 
{\it volume bundle}\index{volume bundle} $(\on{Vol}(M),\pi_M,M)$ of $M$ is
the $1$-dimensional vector bundle (line bundle) 
which is given by the following cocycle of 
transition functions:
\begin{gather*}
\ps_{\al\be}: U_{\al\be}=U_\al\cap U_\be \to \mathbb R\setminus \{0\} = GL(1,\mathbb R),\\
\ps_{\al\be}(x) = |\det d(u_\be\o u_\al\i)(u_\al(x))|
     =\frac1{|\det d(u_\al\o u_\be\i)(u_\be(x))|}.
\end{gather*}
$\on{Vol(M)}$ is a trivial line bundle over $M$.
But there is no natural trivialization.
There is a natural order on each fiber.
Since $\on{Vol}(M)$ is a natural bundle of order 1 on $M$, there is a natural action of the group 
$\Diff(M)$ on $\on{Vol}(M)$, given by 
$$
\xymatrix{
\on{Vol(M)} \ar[d] \ar[rr]^{\quad|\det(T\ph\i)|\o\ph\quad} & &\on{Vol(M)} \ar[d]
\\
M \ar[rr]^{\ph} & & M
}.
$$
If $M$ is orientable, then $\on{Vol}=\La^{m}T^*M$. If $M$ is not orientable, let $\tilde M$ be the 
orientable double cover of $M$ with its deck-transformation $\ta:\tilde M\to \tilde M$. Then 
$\Ga(\on{Vol}(M))$ is isomorphic to the space $\{\om\in\Om^{m}(\tilde M): \ta^*\om=-\om\}$. These are 
the `formes impaires' of de~Rham. See \cite[13.1]{Mic2008} for this.

Sections of the line bundle $\on{Vol}(M)$ 
are called densities. The space $\Ga(\Vol(M))$ of all smooth sections is a Fr\'echet space in its 
natural topology; see \cite{KM1997}.
For each section $\al$ of $\on{Vol}(M)$ of compact support the integral
$\int_M \al$ is invariantly defined as follows:
Let $(U_\al,u_\al)$ be an atlas on $M$ with associated trivialization 
$\ps_\al:\on{Vol}(M)|_{U_\al}\to \mathbb R$, and 
let $f_\al$ be a partition of unity with 
     $\on{supp}(f_\al)\subset U_\al$. Then we put 
\begin{equation*}
\int_M\mu = \sum_\al \int_{U_\al} f_\al\mu 
:= \sum_\al \int_{u_\al(U_\al)} 
     f_\al(u_\al\i(y)).\ps_\al(\mu(u_\al\i(y)))\,dy.
\end{equation*}
The integral is independent of the choice of the atlas and the partition of unity. 

\subsection*{The Fisher--Rao metric}
Let $M^m$ be a smooth compact manifold without boundary.
We denote by $\on{Dens}_+(M)$ the space of smooth positive densities on $M$, i.e. $\on{Dens}_+(M) = \{ \mu \in \Ga(\on{Vol}(M)) \,:\, \mu(x) > 0\; \forall x \in M\}$. Let $\on{Prob}(M)$ be the subspace of positive densities with integral 1 on $M$. Both spaces are smooth Fr\'echet manifolds, in particular they are open subsets of the affine spaces of all densities and densities of integral 1 respectively. For $\mu \in \on{Dens}_+(M)$ we have
$ T_\mu \on{Dens}_+(M) = \Ga(\on{Vol}(M))$ and for $\mu\in \Prob(M)$ we have 
$$
T_\mu\Prob(M)=\{\al\in \Ga(\on{Vol}(M)): \int_M\al =0\}.
$$
The Fisher--Rao metric is a Riemannian metric on $\on{Prob}(M)$ and is defined as follows:
$$
G^{\operatorname{FR}}_\mu(\al,\be) = \int_M \frac{\al}{\mu}\frac{\be}{\mu}\mu.
$$
This metric is invariant under the associated action of $\Diff(M)$ on $\Prob(M)$, since
$$
\Big((\ph^*)^*G^{\operatorname{FR}}\Big)_\mu(\al,\be) = G^{\operatorname{FR}}_{\ph^*\mu}(\ph^*\al,\ph^*\be) 
= \int_M \Big(\frac{\al}{\mu}\o \ph\Big)\Big(\frac{\be}{\mu}\o \ph\Big)\ph^*\mu
= \int_M \frac{\al}{\mu}\frac{\be}{\mu}\mu\,.
$$

The uniqueness result for the Fisher--Rao metric follows from the following classification of $\on{Diff}(M)$-invariant bilinear forms on $\on{Dens}_+(M)$.

\begin{maintheorem*} Let $M$ be a compact manifold without boundary of dimension $\geq 2$.
Let $G$ be a smooth (equivalently, bounded) bilinear form on $\Dens_+(M)$ which is invariant under the action of 
$\Diff(M)$. Then 
$$
G_\mu(\al,\be)=C_1 \int_M \frac{\al}{\mu}\frac{\be}{\mu}\,\mu + C_2\int_M\al \cdot \int_M\be
$$
for some constants $C_1,C_2$.  
\end{maintheorem*}

To see that this theorem implies the uniqueness of the Fisher--Rao metric, note that if $G$ is a $\on{Diff}(M)$-invariant Riemannian metric on $\on{Prob}(M)$, then we can equivariantly extend it to $\on{Dens}_+(M)$ via
\[
G_\mu(\al, \be) = G_{\mu(M)\i \mu}\left( \al - \mu(M) \int_M \al, \be - \mu(M) \int_M \be\right)\,.
\]

\subsection*{Relations to right-invariant metrics on diffeomorphism groups.} 
Let $\mu_0\in \Prob(M)$ be a fixed smooth positive probability density.
In \cite{KLMP2013} it has been shown, that the 
degenerate, $\dot H^1$-metric $\frac12\int_M\on{div}^{\mu_0}(X).\on{div}^{\mu_0}(X).\mu_0$  on 
$\X(M)$ is invariant under the adjoint action of $\Diff(M,\mu_0)$.
Thus the induced degenerate right invariant metric on 
$\on{Diff}(M)$ descends to a metric on  
$$
\Prob(M)\cong \on{Diff}(M,\mu_0)\backslash \on{Diff}(M)\quad\text{ via }\quad\Diff(M)\ni\ph\mapsto \ph^*\mu_0\in\Prob(M)
$$
which is invariant under the right action of $\Diff(M)$. This metric turns out to be 
the Fisher--Rao metric on $\on{Prob}(M)$.
In \cite{Mod2014}, the $\dot H^1$-metric was extended to a non-degenerate metric on $\on{Diff}(M)$, 
that also descends to the Fisher--Rao metric. A consequence of our uniqueness result is the following:

\begin{corollary*}
Let $\operatorname{dim}(M)\geq 2$. If a weak right-invariant (possibly degenerate) Riemannian 
metric $\tilde G$ on $\on{Diff}(M)$ descends to a metric $G$ on $\on{Prob}(M)$, i.e., the map 
$(\on{Diff}(M), \tilde G) \to (\on{Prob}(M), G)$ is a Riemannian submersion, then $G$ has to be a 
multiple of the Fisher--Rao metric.   
\end{corollary*}

For $M=S^1$ the descending property is much less restrictive, since in this case the group of 
volume preserving diffeomorphism is generated by constant vector fields only. 
Thus any right invariant metric on the homogenous space $\Diff(S^1)/S^1$ descends 
to a $\Diff(S^1)$ invariant metric on $\Prob(S^1)$, e.g., the homogenous Sobolev metric of order $n\geq 1$:
\begin{align*}
G_{\operatorname{Id}}(X,Y)=\sum_{k=1}^{n}\int_{S^1} \partial^k_\theta X. \partial^k_\theta Y d\theta\,.
\end{align*}
For $n=1$ the metric descends to the Fisher--Rao metric and for $n=2$ we obtain a higher order metric.
For the one-dimension situation see also the last Section of this article,
where relations between metrics on $\Dens_+(S^1)$ and $\operatorname{Met}(S^1)$ are discussed.

\subsection*{Proof of the Main Theorem.}
Let us first reduce the case of a non-orientable manifold to orientable manifolds. If $M$ is 
non-orientable, let $\tilde M$ be the orientable double cover and $\ta: \tilde M \to \tilde M$ the 
deck-transformation. We can decompose   
\[
\Om^m(\tilde M) = \{ \ta^\ast \om = -\om\} \oplus \{ \ta^\ast \om = \om\}\,,
\]
and $\Dens_+(M)$ is isomorphic to the first summand. Any bilinear form $G$ on $\Dens_+(M)$ can be extended 
to a bilinear form $\tilde G$ on $\Dens_+(\tilde M)$ and the extension is $\Diff(\tilde M)$-invariant. 
Thus we have reduced the proof to the orientable situation.

From now on we assume that $M$ is orientable. Let us fix a basic probability density $\mu_0$. By 
the Moser trick \cite{Moser65}, see \cite[31.13]{Mic2008} or the proof of  
\cite[43.7]{KM1997} for proofs in the notation used here, there exists for each $\mu\in\Dens_+(M)$ a 
diffeomorphism $\ph_\mu\in \Diff(M)$ with $\ph_\mu^*\mu=\mu(M)\mu_0=:c.\mu_0$ where 
$c=\mu(M)=\int_M\mu>0$.
Then 
$$
\big((\ph_\mu^*)^*G\big)_\mu(\al,\be) = G_{\ph_\mu^*\mu}(\ph_\mu^*\al,\ph_\mu^*\be) =  
G_{c.\mu_0}(\ph_\mu^*\al,\ph_\mu^*\be)\,.
$$
Thus it suffices to show that for any $c>0$ we have
$$
G_{c\mu_0}(\al,\be)= \frac{c_1}{c}.\int_M\frac{\al}{\mu_0}\frac{\be}{\mu_0}\mu_0 + 
C_2\int_M\al\cdot\int_M\be
$$
for some constants $c_1,C_2$. 
Both bilinear forms are still invariant under the action of the group
$\Diff(M,c\mu_0)=\Diff(M,\mu_0)=\{\ps\in \Diff(M): \ps^*\mu_0=\mu_0\}$.
The bilinear form 
$$
T_{\mu_0}\Dens_+(M)\x T_{\mu_0}(M)\Dens_+\ni(\al,\be)\mapsto 
G_{c\mu_0}\Big(\frac{\al}{\mu_0}\mu_0,\frac{\be}{\mu_0}\mu_0\Big)
$$
can be viewed as a bilinear form  
\begin{gather*}
C^\infty(M)\x C^\infty(M)\ni (f,g)\mapsto G_c(f,g)\,.
\end{gather*}
We will consider now the associated bounded mapping
$$
\check G_c:C^\infty(M) \to C^\infty(M)' = \mathcal D'(M)\,.
$$

\noindent\thetag{1}
Since we assume that $M$ is orientable, each density is an $m$-form. The Lie algebra $\X(M,\mu_0)$ 
of $\Diff(M,\mu_0)$ consists of vector fields $X$ with  
$\on{div}^{\mu_0}(X) = 0$, or $d i_X \mu_0 = 0$. 
The mapping $\hat\io_{\mu_0}: \X(M)\to \Om^{m-1}(M)$ given by $X\mapsto i_X\mu_0$ is an isomorphism. The Lie 
subalgebra $\X(M,\mu_0)$ of divergence free vector fields corresponds to the space of closed 
$(m-1)$-forms. 
Denote by $\X_{\text{exact}}(M,\mu_0)$ the space of `exact' divergence free vector fields $X = 
\hat\io_{\mu_0}\i(d\om)$, where $\om$ runs through $\Om^{m-2}(M)$. 

\noindent\thetag{2} {\em If for $f\in C^\infty(M)$ and a connected open set $U\subseteq M$ we have 
$\L_X (f|U)=0$ for all 
        $X\in\X_{\text{exact}}(U,\mu_0)$, then $f|U$ is constant.} 

Since we shall need some details later on, we prove this well-known fact.
Let $x \in U$. For every tangent vector $X_x\in T_xM$ we can find a vector field 
$X\in\X_{\text{exact}}(M,\mu_0)$ such that $X(x)=X_x$;
to see this, choose a chart $(U,u)$ near $x$ such that 
$\mu_0|U=du^1\wedge \dots\wedge du^m$, and choose $g\in C_c^\infty(U)$, such 
that $g=1$ near $x$. Then 
$X:= \hat\io_{\mu_0}\i d(g.u^2.du^3\wedge \dots\wedge du^m)\in \X_{\text{exact}}(M,\mu_0)$ and $X = \p_{u^1}$ 
near $x$. So we can produce a basis for $T_xM$ and even a local frame near $x$. 
Thus $\L_X f|U = 0$ for all $X \in \X_{\text{exact}}(M,\mu_0)$ implies $df = 0$ and hence $f$ is constant.

\noindent\thetag{3} {\em If for a distribution $A\in \mathcal D'(M)$ and a connected open set 
$U\subseteq M$ we have $\L_X A|U = 0$ for all $X\in\X_{\text{exact}}(M,\mu_0)$, then $A|U=C\mu_0|U$ 
for some constant $C$, meaning $\langle A, f \rangle = C \int_M f \mu_0$ for all $f \in C_c^\infty(U)$.}

Because $\langle \L_XA,f \rangle = -\langle  A,\L_Xf\rangle$, the invariance property, $\L_XA|U = 
0$, implies $\langle A, \L_X f \rangle = 0$ for all $f\in C^\infty_c(U)$. 
Clearly, $\int_M(\L_Xf)\mu_0 = 0$. Without loss, let us assume now that $U$ is an open chart, that 
is diffeomorphic to $\mathbb R^m$.  
Let $g\in C^\infty_c(U)$ satisfy $\int_M g\mu_0 = 0$; we will show that $\langle A, g \rangle = 0$. 
Because the integral over $g\mu_0$ is zero, the compact cohomology class 
$[g\mu_0]\in H^m_c(U)\cong \mathbb R$ vanishes; thus there exists 
$\al\in \Om^{m-1}_c(U)\subset \Om^{m-1}(M)$ with    
$d\al=g\mu_0$. 
Since we are working on a coordinate chart, which is diffeomorphic to $\R^m$, we can 
write $\al = \sum_j f_j d\be_j$ with $\be_j \in \Om^{m-2}(U)$ and $f_j \in C_c^\infty(U)$. 
Choose $h\in C^\infty_c(U)$ with $h=1$ on $\bigcup_j\on{supp}(f_j)$, so 
that $\al = \sum_j f_j d(h\be_j)$ and $h\be_j\in\Om^{m-2}(M)$.
In 
particular the vector fields $X_j = \hat\io\i_{\mu_0}d(h\be_j)$ lie in $\X_{\text{exact}}(M,\mu_0)$ 
and we have the identity $\sum_{j} f_j.i_{X_i}\mu_0 =\al$.
 This means 
\begin{align*}
\sum_j (\L_{X_j}f_j)\mu_0 &= \sum_j\L_{X_j}(f_j\mu_0) = \sum_j di_{X_j}(f_j\mu_0) = d\Big(\sum_j 
f_j.i_{X_j}\mu_0\Big) = d\al = g\mu_0
\\
\sum_j \L_{X_j}f_j &=  g\,,
\end{align*}
leading to
\[
\langle A,g \rangle = \sum_j \langle A, \L_{X_j} f_j \rangle = -\sum_j \langle \L_{X_j} A, f_j \rangle = 0\,.
\]
So $\langle  A,g\rangle=0$ for all $g\in C^\infty_c(U)$ with $\int_M g\mu_0= 0$. Finally, choose a 
function $\ph$ with support in $U$ and $\int_M \ph \mu_0 = 1$. Then for any $f \in C_c^\infty(U)$, 
the function defined by $g = f - (\int_M f \mu_0). \ph$ in $C^\infty(M)$ satisfies $\int_M g \mu_0 = 0$ and so  
\[
\langle A,f \rangle = \langle A, g \rangle + \langle A, \ph \rangle \int_{M} f \mu_0 = C \int_{M} f \mu_0\,,
\]
with $C = \langle A, \ph\rangle$. Thus $A|U=C\mu_0|U$  and \thetag{3} is proved. 

\noindent\thetag{4} {\em The operator $\check G_c: C^\infty(M)\to \mathcal D'(M)$
has the following property: If for $f\in C^\infty(M)$ and a connected open $U\subseteq M$ the 
restriction $f|U$ is constant, then we have $\check G(f)|U = C_U(f)\mu_0|U$ for some constant $C_U(f)$.}

To see \thetag{4}, for $x\in U$,
choose a smooth function $g$ on $M$ with $g=1$ in a neighborhood of
$M\setminus U$ and $g=0$ on an open neighborhood $V$ of $x$. 
Then for any $X \in \X_{\text{exact}}(M,\mu_0)$, that is $X=\hat\io_{\mu_0}\i(d\om)$ for some 
$\om \in \Om^{m-2}(M)$, let $Y = \hat\io_{\mu_0}\i(d(g\om))$. The vector field $Y$ is again 
divergence free, equals $X$ on a neighborhood of $M\setminus U$, and vanishes on $V$. Since $f$ is 
constant on $U$, it follows that $\L_X f = \L_{Y} f$. Using the invariance of $G_c$, we have for 
all $h \in C^\infty(M)$,   
\[
\left\langle \L_X \check G_c(f), h \right\rangle 
= \left\langle \check G_c(f), -\L_X h \right\rangle
= -G_c(f, \L_X h) = G_c(\L_X f, h)
= \left\langle \check G_c(\L_X f), h \right\rangle\,,
\]
and thus also
\[
\L_X \check G_c(f) = \check G_c(\L_X f) = \check G_c(\L_{Y} f) =
\L_{Y}\check G_c(f)\,.
\]
Now $Y$ vanishes on $V$ and therefore so does $\L_X \check G_c(f)$. By \thetag{3} we have 
$\check G_c(f)|V=C_V(f)\mu_0|V$ for some 
$C_V(f)\in \mathbb R$. Since $U$ is connected, all the constants $C_V(f)$ have to agree, giving a constant $C_U(f)$, depending only on $U$ and $f$.
Thus \thetag{4} follows.

By the Schwartz kernel theorem, $\check G_c$ has a kernel $\hat G_c$, which is
a distribution
(generalized function) in 
\[
\mathcal D'(M\x M) \cong \mathcal D'(M) \bar\otimes \mathcal D'(M) = 
(C^\infty(M)\bar\otimes C^\infty(M))' \cong L(C^\infty(M), \mathcal D'(M))\,.
\] 
Note the defining relations
$$
G_c(f,g) = \langle \check G_c(f),g \rangle = \langle \hat G_c, f\otimes g \rangle.
$$
Moreover, $\hat G_c$ is invariant under the diagonal action of $\Diff(M,\mu_0)$ on $M\x M$. In view 
of the tensor product in the defining relations, the 
infinitesimal version of this invariance is:  $\L_{X\x 0 + 0\x X}\hat G_c = 0$ for all $X\in\X(M,\mu_0)$.

\noindent\thetag{5} {\em 
There exists a constant $C_2$ such that 
the distribution $\hat G_c - C_2 \mu_0\otimes\mu_0$ is supported on the diagonal of $M\x M$.}

Namely, if $(x,y)\in M\x M$ is not on the diagonal, then there exist open neighborhoods $U_x$ of $x$ and 
$U_y$ of $y$ in $M$ such that $\overline{U_x}\x \overline{U_y}$ is disjoint to the diagonal, or 
$\overline{U_x}\cap \overline{U_y}=\emptyset$.  
Choose any functions 
$f,g\in C^\infty(M)$ with $\on{supp}(f)\subset U_x$ and $\on{supp}(g)\subset U_y$. Then 
$f|(M\setminus \overline{U_x}) = 0$, so by \thetag{4}, $\check G_c(f)|(M\setminus \overline{U_x}) = 
C_{M\setminus \overline{U_x}}(f).\mu_0$. Therefore, 
\begin{align*}
G_c(f,g) &= \langle \hat G_c, f\otimes g \rangle = \langle \check G_c(f), g \rangle 
\\&
= \langle \check G_c(f)|(M\setminus \overline{U_x}), g|(M\setminus \overline{U_x})\rangle\,,
&\quad&\text{  since }\on{supp}(g)\subset U_y \subset M\setminus \overline{U_x},
\\&
= C_{M\setminus \overline{U_x}}(f)\cdot\int_M g\mu_0
\end{align*}
By applying the argument for the transposed bilinear form $G_c^T(g,f) = G_c(f,g)$, which is also $\on{Diff}(M,\mu_0)$-invariant, we arrive at
\begin{align*}
G_c(f,g) &= G_c^T(g,f) = C'_{M\setminus \overline{U_y}}(g)\cdot\int_M f\mu_0\,.
\end{align*}
Fix two functions $f_0, g_0$ with the same properties as $f,g$ and additionally $\int_M f_0\mu_0 = 
1$ and $\int_M g_0\mu_0 = 1$. Then we get 
$
C_{M\setminus \overline{U_x}}(f) = C'_{M\setminus \overline{U_y}}(g_0) \int_M f \mu_0\,,
$
and so
\begin{align*}
G_c(f,g) &= C'_{M\setminus \overline{U_y}}(g_0) \int_M f \mu_0 \cdot \int_M g \mu_0
\\&
= C_{M\setminus \overline{U_x}}(f_0) \int_M f \mu_0 \cdot \int_M g \mu_0\,.
\end{align*}
Since $\dim(M)\ge2$ and $M$ is connected, the complement of the diagonal in $M\x M$ is also 
connected, and thus the constants $C_{M\setminus \overline{U_x}}(f_0)$ and 
$C'_{M\setminus \overline{U_y}}(g_0)$ cannot depend on the functions $f_0, g_0$ or the open sets 
$U_x$ and $U_y$ as long as the latter are disjoint. Thus there exists a constant $C_2$ such that for all 
$f,g \in C^\infty(M)$ with disjoint supports we have
$$
G_c(f,g) = C_2\int_M f\mu_0\cdot \int_M g\mu_0
$$
Since $C^\infty_c(U_x\x U_y)=C^\infty_c(U_x)\bar\otimes C^\infty_c(U_y)$, 
this implies claim \thetag{5}. 

Now we can finish the proof. We may replace $\hat G_c\in \mathcal D'(M\x M)$ by $\hat G_c - C_2 
\mu_0\otimes\mu_0$ and thus assume without loss that the constant $C_2$ in \thetag{5} is 0.
Let $(U,u)$ be a chart on $M$ such that $\mu_0|U = du^1\wedge \dots\wedge du^m$.
The distribution $\hat G_c|{U\x U}\in \mathcal D'(U\x U)$ has support contained in the diagonal 
and is of finite order $k$. By \cite[Theorem 5.2.3]{Hor1983}, the corresponding operator 
$\check G_c: C^\infty_c(U)\to \mathcal D'(U)$ is of the form 
$\hat G_c(f) = \sum_{|\al|\le k} A_\al.\p^\al f$ for $A_\al\in \mathcal D'(U)$, so that 
$G(f,g) = \langle \check G_c(f),g \rangle = \sum_{\al} \langle  A_\al,(\p^\al f).g\rangle$. 
Moreover, the $A_\al$ in this representation are uniquely given, 
as is seen by a look at \cite[Theorem 2.3.5]{Hor1983}.

For $x\in U$ choose an open set $U_x$ with $x\in U_x\subset \overline{U_x}\subset U$, and choose 
$X\in \X_{\text{exact}}(M,\mu_0)$ with $X|{U_x}=\p_{u^i}$, as in the proof of \thetag{2}.
For functions $f,g\in C^\infty_c(U_x)$ we then have, by the invariance of $G_c$, 
\begin{align*}
0 &= G_c(\L_Xf,g) + G_c(f,\L_Xg) = \langle \hat G_c|{U\x U}, \L_Xf\otimes g + f\otimes \L_Xg \rangle
\\&
=\sum_\al \langle  A_\al, (\p^\al\p_{u^i} f).g + (\p^\al f)(\p_{u^i}g) \rangle
\\&
=\sum_\al \langle  A_\al, \p_{u^i}((\p^\al f).g) \rangle
=\sum_\al \langle -\p_{u^i} A_\al, (\p^\al f).g \rangle\,.
\end{align*}
Since the corresponding operator has again a kernel distribution which is supported on the 
diagonal, and since the distributions in the representation are unique, we can conclude that 
$\p_{u^i}A_\al|U_x = 0$ for each $\al$, and each $i$.

To see that this implies that $A_\al|U_x = C_{\al} \mu_0|U_x$, let $f \in C_c^\infty(U_x)$ with $\int_M f \mu_0 = 0$. Then, as in \thetag{3}, there exists $\om \in \Om_c^{m-1}(U_x)$ with $d\om = f\mu_0$. 
In coordinates we have $\om = \sum_i \om_i. du^1 \wedge \dots \wedge \widehat{du^i} \wedge du^m$, and so $f = \sum_i (-1)^{i+1} \p_{u^i} \om_i$ with $\om_i \in C_c^\infty(U_x)$. Thus 
\[
\langle A_\al, f \rangle = \sum_i (-1)^{i+1} \langle A_\al, \p_{u^i}\om_i \rangle
= \sum_i (-1)^{i} \langle \p_{u^i} A_\al, \om_i \rangle = 0\,.
\]
Hence $\langle A_\al, f \rangle = 0$ for all $f \in C_c^\infty(U_x)$ with zero integral and as in the proof of \thetag{3} we can conclude that $A_\al|U_x = C_{\al} \mu_0|U_x$.

But then 
$G_c(f,g) = \int_{U_x} (Lf).g \mu_0$
for the differential operator $L = \sum_{|\al|\leq k} C_\al \p^\al$ with constant coefficients on $U_x$.
Now we choose $g\in C_c^\infty(U_x)$ such that $g=1$ on the support of $f$.
By the invariance of $G_c$ we have again
$0= G_c(\L_Xf,g)+G_c(f,\L_Xg)= \int_{U_x} L(\L_X f).\mu_0$ for each $X\in\X(M,\mu_0)$.
Thus the distribution $f\mapsto \int_{U_x} L(f)\mu_0$ vanishes on all functions of the form $\L_Xf$,
and by \thetag{3}
we conclude that $L(\quad).\mu_0= C_x.\mu_0$ in $\mathcal D'(U_x)$, or $L=C_x\on{Id}$.
By covering $M$ with open sets $U_x$, we see that all the constants $C_x$ are the same. 
This concludes the proof of the Main Theorem. \qed

\subsection*{Invariant metrics on $\on{Dens_+}(S^1)$.}
It is interesting to consider the case $M=S^1$, which is not covered by the theorem.
In the following let $M=S^1$. Then positive densities can be represented by positive one-forms. The 
space of all positive densities is isomorphic to the space of all Riemannian metrics on $S^1$ via 
the $\Diff(S^1)$-equivariant mapping  
\begin{align*}
\Phi= (\quad)^2: \on{Dens}_+(S^1) \rightarrow \operatorname{Met}(S^1), 
\qquad 
\Ph(f d\th) = f^2 d\th^2\,.
\end{align*}
On $\operatorname{Met}(S^1)$ there exists a variety of 
$\operatorname{Diff}(S^1)$-invariant metrics; see \cite{BHM2013}. 
We can take for example the family of Sobolev-type metrics. Write $g \in \on{Met}(S^1)$ in the form 
$g = \tilde g d\th^2$ and $h=\tilde h d\th^2$, $k=\tilde k d\th^2$ with $\tilde g,\tilde h, \tilde k\in C^{\infty}(S^1)$. Then for any integer $n$, the following metrics 
are $\on{Diff}(S^1)$-invariant,  
\begin{align*}
 G_{g}^l(h,k) = \int_{S^1} \frac{\tilde h}{\tilde g}.\left(1+\Delta^g\right)^{n}
\left(\frac{\tilde k}{\tilde g}\right) \sqrt{\tilde g} \,d\th\,;
\end{align*}
here $\De^g$ denotes the Laplacian on $S^1$ with respect to the metric $g$. 
Due to the equivariance of $\Phi$, the pullback via $\Phi$ of any of these metrics yields a $\Diff(S^1)$-invariant metric on 
$\Dens_+(M)$, given by
\begin{align*}
G_{\mu}(\alpha,\beta) = 4 \int_{S^1} \frac{\alpha}{\mu} . \left(1+\Delta^{\Phi(\mu)}\right)^{n}\left(\frac{\beta}{\mu}\right)\mu\,. 
\end{align*}
For $n=0$ we obtain 4 times the Fisher--Rao metric. For $n \geq 1$ we see by the number of derivatives involved in the expression for $G_\mu(\al,\be)$, 
that we obtain different $\on{Diff}(S^1)$-invariant metrics on $\on{Dens}_+(M)$ as well as on $\on{Prob}(S^1)$.

\end{document}